\documentclass[12pt, draft]{article}
\usepackage{amsmath}
\usepackage{amsthm}
\usepackage{amsfonts}
\usepackage{amssymb}
\usepackage[dvips]{graphicx}

\newtheorem{theorem}{Theorem}[section]
\newtheorem{utv*}{Proposition}
\newtheorem{hyp*}{Conjecture}

\newtheorem{corollary}[theorem]{Corollary}

\newtheorem*{th*}{Theorem}

\newcommand{\av}[2]{\langle #1\rangle_{_{\scriptstyle #2}}}
\newcommand{\ave}[1]{\langle #1\rangle}

\def\ili{\int\limits}

\def\R{\mathbb{R}}

\def\vf{\varphi}

\begin{document}
\title{An observation: cut-off of the weight $w$ does not increase the $A_{p_{1}, p_{2}}$-``norm'' of $w$}
\author{Alexander Reznikov, \; Vasiliy Vasyunin, \; Alexander Volberg}
\date{}
\maketitle
\begin{abstract}
We consider weights $w$ and their cut-offs: $w_a(t)=w(t)$ if $w(t)\leqslant a$
and $w_a(t)=a$ if $w(t)\geqslant a$. We consider a generalized $A_p$-``norm''
and prove that the ``norm'' of $w_a$ is not greater then the ``norm'' of $w$.
Our  proof in the case $w\in A_2$ is especially simple.
\end{abstract}

\section{Introduction}
Put $I$ be a cube in $\R^n$ and $p_1>p_2$. For every summable function $\vf$
and any subset $J\subset I$ we denote
$$
\av{\vf}{J}=\frac1{|J|}\ili_{J}\vf (t) dt,
$$
where $|J|$ is Lebesgue measure of $J$. For simplicity, when we take an average
over the whole cube $I$, we'll drop the subindex and write $\ave{\vf}$.

Take a nonnegative function $w$. By H\"{o}lder's inequality we have
\begin{equation}
\label{hol}
\av{w^{p_1}}J^{\frac1{p_1}}\av{w^{p_2}}J^{-\frac1{p_2}}\geqslant1\qquad\forall
J\subset [0,1].
\end{equation}
We would like to consider an upper bound of the left-hand side. Precisely, we define
$$
[w]_{p_1,p_2}=
\sup\left(\av{w^{p_1}}J^{\frac1{p_1}}\av{w^{p_2}}J^{-\frac1{p_2}}\right),
$$
where the supremum is taken over all cubes $J$, $J\subset I$. If
$[w]_{p_1,p_2}<\infty$ then we say that $w\in A_{p_1,p_2}$. Note that if
$p_1=1$, $p_2=1-p'=-\frac1{p-1}$ then we get a famous class $A_p$. In this case
instead of $[.]_{1,1-p'}$ we write $[.]_p$.

Notice that the {\it reverse H\"older inequality} is also included as a particular case of these classes $A_{p_1,p_2}$.

We shall assume that $p_i\ne0,\pm\infty$, but it is clear that everything
remains true in the limit cases. In the case $p=0$ the expression
$\av{w^p}J^{\frac1p}$ has to be replaced by $\exp\av{\log w}J$. It has to be
replaced by $\sup_J w$ in the case $p=+\infty$ and by $\inf_J w$ in the case
$p=-\infty$.

We also point out the $A_2$-case: when $p_1=1$ and $p_2=-1$ we have
$$
[w]_2=\sup_{J\subset I}\left( \av wJ\av{w^{-1}}J\right)
$$

Observe that function $w$ can be unbounded or not separated from zero. However,
for some problems it is convenient to consider only bounded, separated from
zero weights.

For a given $a$, $a>0$, we define
$$
w_a(t)=\begin{cases}w(t), &w(t)\leqslant a \\ a, &w(t)\geqslant a \end{cases}.
$$

It was well known that the following inequality is true:
$$
[w_a]_p\leqslant c\cdot [w]_p
$$
with a constant $c$. The main purpose of this text is to delete this constant
and write $1$ instead.

\section{Main results}
We are going to prove the following general theorem.
\begin{theorem}
\label{theorem}
Let $p_1>p_2$; let $w$ be a nonnegative function, defined on $I\subset\R^d$.
Take
$$
w_a(t)=\begin{cases}w(t), &w(t)\leqslant a \\ a, &w(t)> a \end{cases}.
$$
Then for every cube $J$, $J\subset I$, the following is true:
\begin{equation}\label{mainineq}
\av{w_a^{p_1}}J^{\frac1{p_1}}\av{w_a^{p_2}}J^{-\frac1{p_2}}-
\av{w^{p_1}}J^{\frac1{p_1}}\av{w^{p_2}}J^{-\frac1{p_2}}\leqslant 0.
\end{equation}
Consequently,
$$
[w_a]_{p_1,p_2}\leqslant [w]_{p_1,p_2}.
$$
\end{theorem}
This theorem gives an answer to a similar question, when we cut from below. Precisely,
\begin{corollary}\label{coll1}
Denote
$$
w^a(t)=\begin{cases}w(t), &w(t)\geqslant a \\ a, &w(t)<a \end{cases}.
$$
Then the following inequality holds:
$$
[w^a]_{p_1,p_2}\leqslant [w]_{p_1,p_2}.
$$
\end{corollary}
This corollary is an immediate consequence of the theorem, since instead of
$a,w,p_1,p_2$ we can consider $\frac1a,\frac1w,-p_2,-p_1$.

\begin{corollary}\label{coll}
Take a function $w\in A_{p_{1}, p_{2}}$.
For every integer $n$, $n\geqslant 1$, denote
$$
\vf_{n}(t)=
\begin{cases}
n, &w(t)>n \\
w(t), &\frac1n<w(t)\leqslant n \\
\frac1n, &w(t)\leqslant n
\end{cases}.
$$
Then
\begin{align}
[\vf_n]_{p_1,p_2}\leqslant [w]_{p_1,p_2}\label{one} \\
\lim_{n\to \infty} [\vf_n]_{p_1,p_2}=[w]_{p_1,p_2}\label{two}.
\end{align}
\end{corollary}

We give an independent proof of~\eqref{mainineq} in one leading particular cases of the class $A_2$. 

The possibility to approximate a function in the class $A_p$ by bounded functions from the same class and with the control of their $A_p$ constants (and we have the best possible control here) 
can be used in various places. First of all, \cite{BM} shows how this can be used to show that the set $\{p: w\in A_p\}$ is open. Secondly, working with Bellman function proofs of various sharp
reverse H\"older inequalities or sharp John--Nirenberg type inequalities (see e.g. \cite{R}, \cite{V}), one needs an approximation of a weight $w$ in $A_p$ (and, more generally, $w\in A_{p_1,p_2}$) by the weights bounded from above and from below and of at most the same $A_p$ constant. We show how this can be easily achieved by just a standard ``cut-off" procedure on weights. Seems like this has not been observed in the literature, even though it amounts to a very simple remark.

\section{History of the question}
As far as we know, the known result for $w_a$, $w^a$ and $\vf_n$ is the
following inequality:
$$
[w_a]_{p_1,p_2}\leqslant 2 [w]_{p_1,p_2}.
$$
In this work we erase the constant $2$. We should cite the work ~\cite{BM},
where the different approach is described. Authors consider weights
$$
\frac{s+w(t)}{s^{2}+sw(t)+1},
$$
which are bounded and which also satisfy~\eqref{one} and~\eqref{two} as $s\to
+0$. However, we think that our approach is more natural if one wants to get a
bounded weight separated from zero.

\section{Proof of the Corollary \ref{coll}}
Inequality~\eqref{one} follows from the main theorem~\ref{theorem} and from the
corollary~\ref{coll1}. Thus, we need to prove~\eqref{two}. By the monotone
convergence theorem we have
$$
\av{\vf_n^{p_k}}J \to \av{w^{p_k}}J,
$$
therefore, for every $J\subset I$, the following is true:
$$
\av{\vf_n^{p_1}}J^{\frac1{p_1}}\av{\vf_n^{p_2}}J^{-\frac1{p_2}}\to
\av{w^{p_1}}J^{\frac1{p_1}}\av{w^{p_2}}J^{-\frac1{p_2}}, \quad n\to\infty.
$$
Therefore,
\begin{align}
\av{w^{p_1}}J^{\frac1{p_1}}\av{w^{p_2}}J^{-\frac1{p_2}}&=
\lim\left(\av{\vf_n^{p_1}}J^{\frac1{p_1}}\av{\vf_n^{p_2}}J^{-\frac1{p_2}}\right)=
\liminf\left(\av{\vf_n^{p_1}}J^{\frac1{p_1}}\av{\vf_n^{p_2}}J^{-\frac1{p_2}}\right)
\\
&\leqslant\liminf
[\vf_n]_{p_1,p_2}\leqslant\limsup[\vf_n]_{p_1,p_2}\leqslant[w]_{p_1,p_2}.
\end{align}
Passing to the supremum over $J$ in the right-hand side, we get
$$
[w]_{p_1,p_2}\leqslant\liminf[\vf_n]_{p_1,p_2}\leqslant
\limsup[\vf_n]_{p_1,p_2}\leqslant[w]_{p_1,p_2},
$$
which finishes the proof.

In next three sections we prove the main Theorem~\ref{theorem}.

\section{The case of $A_2$-weights}
We separate this case since here everything is in some sense linear, and we can
prove everything without taking derivatives. In this case $p_1=1$, $p_2=-1$.
Fix a cube $J\subset I$ and denote
$$
J_1=\{t\in J\colon w(t)\leqslant a\},\qquad J_2=\{t\in J\colon w(t)>a\},
$$
$$
x_i=\av w{J_i},\qquad y_i=\av{\frac1w}{J_i},\qquad\alpha_i=\frac{|J_i|}{|J|}.
$$
Then
\begin{align*}
\av wJ\av{w^{-1}}J&-\av{w_a}J\av{w_a^{-1}}J
\\
&=(\alpha_1x_1+\alpha_2x_2)(\alpha_1y_1+\alpha_2y_2)-
(\alpha_1x_1+\alpha_2a)(\alpha_1y_1+\alpha_2a^{-1})
\\
&=\alpha_1\alpha_2(x_1y_2+x_2y_1-y_1a-x_1a^{-1})+\alpha_2^2(x_2y_2-1).
\end{align*}
The expression in the second parentheses is positive and therefore it is
sufficient to check that the expression in the first parentheses is positive as
well.
\begin{align*}
x_1y_2&+x_2y_1-y_1a-x_1a^{-1}=x_1(y_2-a^{-1})+y_1(x_2-a)
\\
&=\av{x_1(w^{-1}-a^{-1})+y_1(w-a)}{J_2}= \av{\frac{w-a}{wa}(way_1-x_1)}{J_2}.
\end{align*}
Since $y_1\geqslant\frac1a$ and $x_1\leqslant a$, we have $way_1-x_1\geqslant w-a$, which
finishes the proof.

\section{Proof of the general case}
In this section we present a fully general proof.

We keep all notation from the preceding section with a natural modification. Fix a cube $J\subset I$ and put
$$
J_1=\{t\in J\colon w(t)\leqslant a\},\qquad J_2=\{t\in J\colon w(t)>a\},
$$
$$
x_i=\av{w^{p_1}}{J_i},\qquad y_i=\av{w^{p_2}}{J_i},\qquad\alpha_i=\frac{|J_i|}{|J|}\,.
$$
Then we want to prove
\begin{equation}
\label{phi}
\begin{aligned}
\av{&w^{p_1}}J^{\frac1{p_1}}\av{w^{p_2}}J^{-\frac1{p_2}}- \av{w_a^{p_1}}J^{\frac1{p_1}}\av{w_a^{p_2}}J^{-\frac1{p_2}}
\\
&=(\alpha_1x_1+\alpha_2x_2)^{\frac1{p_1}} (\alpha_1y_1+\alpha_2y_2)^{-\frac1{p_2}}-
(\alpha_1x_1+\alpha_2a^{p_1})^{\frac1{p_1}} (\alpha_1y_1+\alpha_2a^{p_2})^{-\frac1{p_2}}
\geqslant0\,.
\end{aligned}
\end{equation}

By H\"older's inequality, we get $x_i^{\frac1{p_1}}\geqslant y_i^{\frac1{p_2}}$. Therefore, if we denote $y_2^{\frac1{p_2}}$ by $u$, then $x_2^{\frac1{p_1}}=su$ for a number $s\geqslant 1$ and expression~\eqref{phi} we need to estimate can be written as the following function of $s$ and $u$:
$$
\vf(s,u)=(\alpha_1x_1+\alpha_2s^{p_1}u^{p_1})^{\frac1{p_1}} (\alpha_1y_1+\alpha_2u^{p_2})^{-\frac1{p_2}}-
(\alpha_1x_1+\alpha_2a^{p_1})^{\frac1{p_1}} (\alpha_1y_1+\alpha_2a^{p_2})^{-\frac1{p_2}}\,.
$$
Since
$$
\frac{\partial\vf}{\partial s}=\alpha_2 s^{p_1-1}u^{p_1} (\alpha_1x_1+\alpha_2s^{p_1}u^{p_1})^{\frac1{p_1}-1}\geqslant0\,,
$$
the function $\vf$ is increasing in $s$ and therefore $\vf(s,u)\geqslant\vf(1,u)$, i.e., it has the minimal value when $w(t)$ is equal to $u$ on $J_2$ identically.

Now we have $u=w(t)|_{J_2}>a$ and since $\vf(1,a)=0$, the desired inequality will be proved after checking that $\frac{\partial\vf}{\partial u}(1,u)\geqslant0$.

\begin{align*}
\frac{\partial\vf}{\partial u}&(1,u)
\\
&=\alpha_2u^{-1}(\alpha_1x_1+\alpha_2u^{p_1})^{\frac1{p_1}-1} (\alpha_1y_1+\alpha_2u^{p_2})^{-\frac1{p_2}-1}\times
\\
&\qquad\qquad\qquad\qquad\qquad\qquad\qquad \times\big[u^{p_1}(\alpha_1y_1+\alpha_2u^{p_2})- u^{p_2}(\alpha_1x_1+\alpha_2u^{p_1})\big]
\\
&=\alpha_1\alpha_2u^{-1}(\alpha_1x_1+\alpha_2u^{p_1})^{\frac1{p_1}-1} (\alpha_1y_1+\alpha_2u^{p_2})^{-\frac1{p_2}-1}[u^{p_1}y_1-u^{p_2}x_1]
\end{align*}
and we are done because $u^{p_1}y_1-u^{p_2}x_1\geqslant0$. Indeed,
since $u\geqslant w(t)$ and $p_1\geqslant p_2$, we have $u^{p_1-p_2}\geqslant w(t)^{p_1-p_2}$, whence $u^{p_1}w^{p_2}\geqslant u^{p_2}w^{p_1}$. Therefore,
$$
u^{p_1}y_1-u^{p_2}x_1=\av{u^{p_1}w^{p_2}-u^{p_2}w^{p_1}}{J_1}\geqslant0\,,
$$
what completes the proof.

\end{document}